\documentstyle[amssymb,amsfonts]{amsart}

\def\C{{\hbox{\bf C}}}

\font \roman = cmr10 at 10 true pt

\def\bas{\begin{align*}}
\def\eas{\end{align*}}
\def\bi{\begin{itemize}}
\def\ei{\end{itemize}}

\def\supp{{\hbox{\roman supp}}}

\def\Z{{\hbox{\bf Z}}}

\newenvironment{proof}{\noindent {\bf Proof} }{\endprf\par}
\def \endprf{\hfill  {\vrule height6pt width6pt depth0pt}\medskip}
\def\emph#1{{\it #1}}
\def\textbf#1{{\bf #1}}

\theoremstyle{plain}
  \newtheorem{theorem}[subsection]{Theorem}

  \newtheorem{lemma}[subsection]{Lemma}
  \newtheorem{corollary}[subsection]{Corollary}

\theoremstyle{remark}

\theoremstyle{definition}

\include{psfig}

\begin{document}

\title[Uncertainty principle for $\Z/p\Z$]{An uncertainty principle for cyclic groups of prime order}

\author{Terence Tao}
\address{Department of Mathematics, UCLA, Los Angeles CA 90095-1555}
\email{tao@@math.ucla.edu}


\begin{abstract}  Let $G$ be a finite abelian group, and let $f: G \to \C$ be a complex
function on $G$.  The uncertainty principle asserts that the support $\supp(f) := \{ x \in G: f(x) \neq 0\}$
is related to the support of the Fourier transform $\hat f: G \to \C$ by the formula
$$ |\supp(f)| |\supp(\hat f)| \geq |G|$$
where $|X|$ denotes the cardinality of $X$.  In this note we show that when $G$ is the cyclic group $\Z/p\Z$
of prime order $p$, then we may improve this to
$$ |\supp(f)| + |\supp(\hat f)| \geq p+1$$
and show that this is absolutely sharp.  As one consequence, we see that a sparse polynomial in $\Z/p\Z$ consisting
of $k+1$ monomials can have at most $k$ zeroes.  Another consequence is a short proof of the 
well-known Cauchy-Davenport inequality.
\end{abstract}

\maketitle

\section{Introduction}

Let $G$ be a finite abelian additive group, and let $e: G \times G \to S^1 := \{ z \in \C: |z| = 1\}$ be any 
non-degenerate \emph{bi-character} of $G$, by which we mean a function $e(x,\xi)$ taking values on the unit circle
obeying the multiplicativity properties
$$ e(x+x',\xi) = e(x,\xi) e(x',\xi); \quad e(x,\xi+\xi') = e(x,\xi) e(x,\xi')$$
and is non-degenerate in the sense that for every $x \neq 0$ there exists a $\xi \in G$ such that $e(x,\xi) \neq 1$,
and conversely for every $\xi \neq 0$ there exists an $x \in G$ such that $e(x,\xi) \neq 1$.  For instance,
if $G$ is the cyclic group $G := \Z/N\Z$, we may take $e(x,\xi) := e^{2\pi i x\xi/N}$.  If $f: G \to \C$ is any
complex-valued function on $G$, we may then define the Fourier transform $\hat f: G \to \C$ by the formula
$$ \hat f(\xi) := \frac{1}{|G|} \sum_{x \in G} f(x) \overline{e(x,\xi)},$$
where $|G|$ denotes the cardinality of $G$.  If we use $\supp(f) := \{ x \in G: f(x) \neq 0\}$ to denote
the support of $f$, we thus see from the triangle inequality, Cauchy-Schwarz and Plancherel that
\begin{align*}
\sup_{\xi \in G} |\hat f(\xi)| 
&\leq \frac{1}{|G|} \sum_{x \in G} |f(x)| \\
&\leq \frac{|\supp(f)|^{1/2}}{|G|^{1/2}} (\frac{1}{|G|} \sum_{x \in G} |f(x)|^2)^{1/2} \\
&= \frac{|\supp(f)|^{1/2}}{|G|^{1/2}} (\sum_{\xi \in G} |\hat f(\xi)|^2)^{1/2} \\
&\leq \frac{|\supp(f)|^{1/2} |\supp(\hat f)|^{1/2}}{|G|^{1/2}} \sup_{\xi \in G} |\hat f(\xi)|.
\end{align*}
Thus, if $f$ is non-zero, we thus obtain the well-known \emph{uncertainty principle}\cite{donoho}, \cite{smith}
\begin{equation}\label{basic-uncertainty} 
|\supp(f)| |\supp(\hat f)| \geq |G|.
\end{equation}
This bound is of course sharp when $f$ is a Dirac mass, or when $\hat f$ is a Dirac mass.  More generally,
if $H$ is any subgroup of $G$, and we set $f$ to be the characteristic function $\chi_H$ of $f$, it is
easy to see that $|\supp(f)| = |H|$ and $|\supp(\hat f)| = |G|/|H|$, so again \eqref{basic-uncertainty} 
is sharp.  Indeed one can show that up to the symmetries of the Fourier transform (translation, modulation,
and homogeneity) this is the only way in which \eqref{basic-uncertainty} can be obeyed with equality (see e.g. 
\cite{prez}).  For more background on the Fourier transform on finite abelian groups and the uncertainty 
principle we refer to \cite{terras}.

Now consider the case where $G$ is a cyclic group of prime order, $G = \Z/p\Z$, with $e(x,\xi) := e^{2\pi i x \xi/p}$.
In this case $G$ has no subgroups other than the trivial ones $\{0\}$ and $G$, and thus one expects to
improve upon \eqref{basic-uncertainty}.  Indeed we can get an absolutely sharp result as to the possible values of
$\supp(f)$ and $\supp(\hat f)$:

\begin{theorem}\label{uncertainty}  Let $p$ be a prime number.
  If $f: \Z/p\Z \to \C$ is a non-zero function, then\footnote{This inequality was also discovered 
independently by Andr\'as Bir\'o \cite{frenkel} and Roy Meshulam (Vsevolod Lev, private
communication).  Given the number of times Lemma \ref{minors} appears to have been rediscovered in the literature it is in fact quite likely that this theorem has existed previously in some unpublished form.}
$$ |\supp(f)| + |\supp(\hat f)| \geq p+1.$$
Conversely, if $A$ and $B$ are two non-empty subsets of $\Z/p\Z$ such that $|A| + |B| \geq p+1$, then there exists
a function $f$ such that $\supp(f) = A$ and $\supp(\hat f) = B$.
\end{theorem}

The informal explanation of this principle is that the class of functions $f$ from $\Z/p\Z \to \C$ has exactly
$p$ degrees of freedom.  Requiring that $\supp(f) = A$ takes away $p-|A|$ of these degrees, while
requiring that $\supp(\hat f) = B$ takes away another $p-|B|$.  The uncertainty principle is then a statement
that the Fourier basis (of exponentials) and the physical space basis (of Dirac deltas) are ``totally skew''
(or more precisely, that all the minors of the exponential basis matrix $(e^{2\pi i jk/p})_{0 \leq j,k < p}$ are
non-zero).  The idea that the prime cyclic group $\Z/p\Z$ has this ``maximally skew'' structure (in some sense,
it is as far as possible from containing subgroups) is consistent with some other recent work on
the arithmetic structure of prime cyclic groups, see for instance \cite{bkt}, \cite{k}.

The proof of Theorem \ref{uncertainty} requires a number of preliminary lemmas.  We first need a lemma from the 
Galois theory of the cyclotomic integers.

\begin{lemma}\label{galois}  Let $p$ be a prime, $n$ be a positive integer, and let
$P(z_1, \ldots, z_n)$ be a polynomial with integer co-efficients.  Suppose that we have $n$ $p^{th}$ roots
of unity $\omega_1, \ldots, \omega_n$ (not necessarily distinct) such that $P(\omega_1, \ldots, \omega_n) = 0$.
Then $P(1, \ldots, 1)$ is a multiple of $p$.
\end{lemma}

\begin{proof}  Write $\omega := e^{2\pi i/p}$, then for every $1 \leq j \leq n$ we
have $\omega_j = \omega^{k_j}$ for some integers $0 \leq k_j < p$.  If we then define the single-variable
polynomial $Q(z)$ by
$$ Q(z) := P(z^{k_1}, \ldots, z^{k_n}) \mod z^p-1,$$
where $R(z) \mod z^p - 1$ is the remainder when dividing $R(z)$ by $z^p - 1$ (or equivalently, taking
the polynomial $R(z)$ and replacing $z^{qp+r}$ with $z^r$ for all $q \geq 1$ and $0 \leq r < p$), then
we have $Q(\omega) = 0$ and $Q(1) = P(1, \ldots, 1)$.  But $Q(z)$ is a polynomial of degree at most $p-1$
with integer coefficients,
and thus must be an integer multiple of the minimal polynomial $1 + z + \ldots + z^{p-1}$ of $\omega$.
The claim follows.
\end{proof}

Using this lemma, we can show that all the minors of the Fourier matrix are non-zero.

\begin{lemma}\label{minors}  Let $p$ be a prime and $1 \leq n \leq p$.  Let $x_1, \ldots, x_n$ be distinct elements
of $\Z/p\Z$, and let $\xi_1, \ldots, \xi_n$ also be distinct elements of $\Z/p\Z$.  Then the matrix
$( e^{2\pi i x_j \xi_k/p} )_{1 \leq j,k \leq n}$ has non-zero determinant.
\end{lemma}

This result was first proved by Chebotar\"ev in 1926 (see \cite{sl}), and with additional proofs given by 
Resetnyak \cite{res}, Dieudonn\'e \cite{d}, Newman \cite{newman}, Evans and Stark \cite{ei}, and more recently
Frenkel \cite{frenkel} and Goldstein, Guralnick, and Isaacs \cite{ggi}.  Recently, some more quantitative measure of the non-degeneracy of (randomly selected) minors was obtained in \cite{candes}.  All proofs of Lemma \ref{minors} 
require a certain amount of algebraic information about the cyclotomic integers, but our proof requires 
relatively little in that regard (all we need is Lemma \ref{galois}).

\begin{proof}
Write $\omega_j := e^{2\pi i x_j/p}$.  Then each $\omega_j$ is a distinct root of
unity, and our task is to show that
$$ \det ( \omega_j^{\xi_k} )_{1 \leq j,k \leq n}$$
is non-zero.  Motivated by the previous lemma, we define the polynomial $D(z_1, \ldots, z_n)$ of 
$n$ variables by
$$ D(z_1, \ldots, z_n) := \det ( z_j^{\xi_k} )_{1 \leq j,k \leq n};$$
here we identify the frequencies $\xi_k \in \Z/p\Z$ with elements of $\{0,1,\ldots,p-1\}$ in
the obvious manner.  This is clearly a polynomial with integer co-efficients.   Unfortunately
$D(1,\ldots,1)$ degenerates to zero and so Lemma \ref{galois} does not
directly tell us that $D(\omega_1, \ldots,\omega_n)$ is non-zero.  Indeed,
$D$ clearly vanishes when $z_j = z_{j'}$ for any $1 \leq j < j' \leq n$, and
so we can factor
\begin{equation}\label{d-factor}
 D(z_1, \ldots, z_n) = P(z_1, \ldots, z_n) \prod_{1 \leq j < j' \leq n} (z_j - z_{j'})
\end{equation}
for some other polynomial $P$ with integer coefficients.  We will show that
$P(1,\ldots,1)$ is not a multiple of $p$, which by Lemma \ref{galois} shows
that $P(\omega_1, \ldots, \omega_n)$ is non-zero, which proves the claim since
the $\omega_j$ are all distinct.

To compute $P(1,\ldots, 1)$, we differentiate $D$ repeatedly.  In particular, we consider the expression
\begin{equation}\label{big-mess}
 (z_1 \frac{d}{dz_1})^0 (z_2 \frac{d}{dz_2})^1 \ldots (z_n \frac{d}{dz_n})^{n-1} D(z_1, \ldots,z_n)|_{z_1 = \ldots = z_n = 1}.
\end{equation}
Note that we are applying $0 + 1 + \ldots + n-1 = \frac{n(n-1)}{2}$ differentiation operators, which
is exactly the same number as the number of linear factors $(z_j - z_{j'})$ in \eqref{d-factor}.  By the Leibnitz rule,
each differentiation operator $z_j \frac{d}{dz_j}$ may eliminate one of these linear factors (and replace it 
with $z_j$), or it may differentiate some other term (e.g. it may differentiate $P$).  But the only terms from
the Leibnitz expansion which do not vanish when $z_1 = \ldots = z_n = 1$ are those in which every
differentiation operator eliminates one of the linear factors (so in particular we never need to
differentiate $P$).  The $n-1$ copies of the differentiation operators
$z_n \frac{d}{dz_n}$ can only eliminate the $n-1$ linear factors $(z_j - z_n)$, and so every one of those linear
factors must be eliminated by one of these differentiation operators, and there are $(n-1)!$ ways in which
this can occur.  We then argue similarly with the $n-2$ copies of $z_{n-1} \frac{d}{dz_{n-1}}$, which
must eliminate the $n-2$ linear factors $(z_j - z_{n-1})$ (and there are $(n-2)!$ ways of doing so).  
Continuing in this fashion we thus see that 
$$ \eqref{big-mess} = (n-1)! (n-2)! \ldots 0! P(1,\ldots,1).$$
Since $(n-1)! (n-2)! \ldots 0!$ is not a multiple of $p$, it thus suffices to show that the integer
\eqref{big-mess} is not a multiple of $p$.  But by the definition of $D(z_1, \ldots, z_n)$ and the
multilinearity of the determinant, and the trivial observation that $(z_j \frac{d}{dz_j}) z_j^\xi = 
\xi z_j^\xi$, we see that
$$ \eqref{big-mess} = \det( \xi_k^{j-1} )_{1 \leq j,k \leq n}.$$
This is a Vandermonde determinant which (as is well-known) is equal to
$$ \pm \prod_{1 \leq k < k' \leq n} (\xi_k - \xi_{k'}).$$
But since the $\xi_k$ are all distinct modulo $p$, this is not a multiple of $p$, and the claim follows.
\end{proof}

From the above Lemma we immediately obtain

\begin{corollary}\label{n-injective}  If $p$ is a prime, and $A, \tilde A$ are non-empty subsets of
$\Z/p\Z$ such that $|A|=|\tilde A|$, then the linear transformation $T: l^2(A) \to l^2(\tilde A)$ defined by
$Tf = \hat f|_{\tilde A}$ (i.e. we restrict the Fourier transform of $f$ to $\tilde A$) is invertible.
Here we use $l^2(A)$ to denote those functions $f: G \to \C$ which are equal to zero outside of $A$.
\end{corollary}

Indeed, the coefficient matrix of $T$ is of the form considered in Lemma \ref{minors}.  From this Corollary
we can now easily prove the uncertainty principle.

\begin{proof}[of Theorem \ref{uncertainty}.]  Suppose for contradiction that we had a non-zero function
$f$ such that $|\supp(f)| + |\supp(\hat f)| \leq p$.  Then if we write $A := \supp(f)$, then
we can find a set $\tilde A$ in $\Z/p\Z$ which is disjoint from $\supp(\hat f)$ and has cardinality equal to $|A|$.  But
this contradicts Corollary \ref{n-injective} since $Tf = 0$ but $f \neq 0$.

Now we prove the converse.  It will suffice to prove the claim when $|A| + |B| = p+1$, since the claim for
$|A| + |B| > p+1$ then follows by applying the claim to subsets $A'$, $B'$ of $A$, $B$ respectively for which
$|A'| + |B'| = p+1$, and then taking generic linear combinations as $A'$, $B'$ vary.

We can then choose an $\tilde A$ in $\Z/p\Z$ of cardinality $|\tilde A|=|A|$ such that $\tilde A$ 
intersects $B$ in only one point,
say $\tilde A \cap B = \{\xi\}$.  But by Corollary \ref{n-injective} the map $T$ is invertible, and in particular
we can find a non-zero $f \in l^2(A)$ such that $\hat f$ vanishes on $\tilde A \backslash \{\xi\}$ and is non-zero
on $\xi$.  Such a function must then be non-zero on all of $A$ and non-zero on all of $B$ since this would
violate the first part of the uncertainty principle proven in the previous paragraph.  Thus $\supp(f) = A$ and
$\supp(\hat f) = B$ as desired.
\end{proof}

Observe that an immediate consequence of Theorem \ref{uncertainty} is that any sparse 
polynomial $\sum_{j=0}^k c_j z^{n_j}$ with $k+1$ non-zero coefficients and $0 \leq n_0 < \ldots < n_k < p$, 
when restricted to the $p^{th}$ roots of unity $\{ z: z^p = 1\}$, can have at most
$k$ zeroes.  Indeed, such a polynomial is essentially the Fourier transform in $\Z/p\Z$ of a function whose 
support has cardinality $k+1$, and so the support of the polynomial must contain at least $p-k$ $p^{th}$ roots
of unity by Theorem \ref{uncertainty}, and the claim follows.  

Another immediate consequence is the Cauchy-Davenport inequality \cite{cauchy}, \cite{davenport}, 
which asserts that for any two finite non-empty
subsets $A$ and $B$ of $\Z/p\Z$, the sumset $A+B := \{a+b: a \in A, b \in B\}$ obeys the bounds
$$ |A+B| \geq \min(|A|+|B|-1, p).$$

\begin{proof}\footnote{We thank Robin Chapman for this proof, which is slightly shorter than the original
proof of the author.}  Fix $A$, $B$.  Since $A$ and $B$ are non-empty, we may find
two subsets $X$ and $Y$ of $\Z/p\Z$ such that $|X| = p+1-|A|$, $|Y| = p+1-|B|$, and $|X \cap Y| = \max(|X| + |Y| - p,1)$.  By Theorem \ref{uncertainty} we may find a function $f$ such that
$\supp(f) = A$ and $\supp(\hat f) = X$, and a function $g$ such that $\supp(g) = B$ and $\supp(\hat g) = Y$.  Then
$f*g$ has support contained in $A+B$ and has Fourier support equal to $X \cap Y$ (in particular, $f*g$ is non-zero), 
and hence by Theorem \ref{uncertainty} again we have $|A+B| + |X \cap Y| \geq p + 1$, which 
gives $|A+B| \geq \max(|A| + |B| - 1,p)$ as desired. 
\end{proof}

It is interesting to compare this proof with the polynomial method proof of \cite{alon}, which uses the basis of
polynomials rather than the basis of exponentials but is otherwise rather similar in spirit.

Based on this result for groups of prime order, it seems natural to conjecture that one can improve
\eqref{basic-uncertainty} substantially for all finite abelian groups $G$, provided that the cardinality of
$|\supp(f)|$ and $|\supp(\hat f)|$ stays well away from any factor of $|G|$.  For instance, Roy Meschulam (private
communication) has used Theorem \ref{uncertainty} and an iteration argument to show that 
$p^j |\supp(f)| + p^{n-j-1} |\supp (\hat f)| \geq p^n + p^{n-1}$ for all non-zero functions $f$ supported
on $(\Z/p\Z)^n$ and all $0 \leq j \leq n-1$.  To put this another way, the point $(|\supp(f)|, |\supp(\hat f)|)$ in
$\Z \times \Z$ lies on or above the convex hull of the points $(p^j, p^{n-j})$ for $0 \leq j \leq n$, which
correspond to the cases where $f$ is the characteristic function of a subgroup of $(\Z/p\Z)^n$.
This has immediate application to the number of zeroes of a sparse
 polynomial of several variables in $\Z/p\Z$, which may be of interest for cryptographic applications.    

\section{Acknowledgements}

This work was conducted at Australian National University.  The author is a Clay Prize Fellow and is supported
by a grant from the Packard Foundation.  The author is also indebted to Robin Chapman for pointing out the
provenance of Lemma \ref{minors} and simplifying the proof of the Cauchy-Davenport inequality, to 
Michael Cowling and Gerd Mockenhaupt for pointing out the provenance of \eqref{basic-uncertainty},
and to Roy Meshulam to pointing out extensions of Theorem \ref{uncertainty} to higher powers of $\Z/p\Z$.
We also thank Gergely Harcos, Melvyn Nathanson and Vselvolod Lev for some corrections and comments.

\end{document}